\documentclass{amsart}
\usepackage{amssymb}
\usepackage{amsmath}

\begin{document}
\def\id{\operatorname{id}}
\def\gronk{\vphantom{\vrule height 12pt}}
\def\mapright#1{\smash{\mathop{\longrightarrow}\limits\sp{#1}}}
\newtheorem{theorem}{Theorem}[section]
\newtheorem{lemma}[theorem]{Lemma}
\newtheorem{remark}[theorem]{Remark}
\newtheorem{definition}[theorem]{Definition}
\newtheorem{corollary}[theorem]{Corollary}
\newtheorem{example}[theorem]{Example}
\def\qedbox{\hbox{$\rlap{$\sqcap$}\sqcup$}}
\makeatletter
  \renewcommand{\theequation}{%
   \thesection.\alph{equation}}
  \@addtoreset{equation}{section}
 \makeatother
\def\BB{\mathcal{B}}
\def\MM{{\mathfrak{M}}}
\title{Jacobi--Tsankov manifolds which are not $2$-step nilpotent}
\author{M. Brozos-V\'azquez}
\address{MBV: Department of Geometry and Topology, Faculty of Mathematics, University of
Santiago de Compostela, 15782 Santiago de Compostela, Spain\vspace{1mm}\\ Email:\it mbrozos@usc.es.}
\author{P. Gilkey}
\address{PG: Mathematics Department, University of Oregon, Eugene, OR 97403, USA\vspace{1mm}\\ Email:
\it gilkey@uoregon.edu.}
\author{S. Nik\v cevi\'c}
\address{SN: Mathematical Institute, SANU,
Knez Mihailova 35, p.p. 367,
11001 Belgrade, 
Serbia and Montenegro\vspace{1mm}\\Email: \it stanan@mi.sanu.ac.yu.}
\begin{abstract} There is a $14$-dimensional algebraic curvature tensor
which is Jacobi--Tsankov (i.e. $\mathcal{J}(x)\mathcal{J}(y)=\mathcal{J}(y)\mathcal{J}(x)$ for all $x,y$) but which is
not
$2$-step Jacobi nilpotent (i.e. $\mathcal{J}(x)\mathcal{J}(y)\ne0$ for some $x,y$); the minimal dimension where this
is possible is $14$. We determine
the group of symmetries of this tensor and show that it is geometrically
realizable by a wide variety of pseudo-Riemannian manifolds which are
geodesically complete and have vanishing scalar invariants. Some of the manifolds in the family are
symmetric spaces. Some are $0$-curvature homogeneous but not locally homogeneous.
\end{abstract}

\maketitle

\section{Introduction}

Let $\nabla$, $\mathcal{R}$, $R$, and $\mathcal{J}$ denote the
Levi-Civita connection, the curvature operator, the curvature tensor, and the Jacobi operator, respectively, of a
pseudo-Riemannian manifold $\mathcal{M}:=(M,g)$:
\begin{eqnarray*}
&&\mathcal{R}(X,Y)=\nabla_X\nabla_Y-\nabla_Y\nabla_X-\nabla_{[X,Y]},\\
&&R(X,Y,Z,W):=g(\mathcal{R}(X,Y)Z,W),\quad\text{and}\\&&\mathcal{J}(x)y:=\mathcal{R}(y,x)x\,.
\end{eqnarray*}

The relationship between algebraic properties of the Jacobi operator and the underlying geometry of the
manifold has been extensively studied in recent years. For example, $\mathcal{M}$ is said to be {\it
Osserman} if the eigenvalues of $\mathcal{J}$ are constant on the pseudo-sphere bundles $S^\pm(M,g)$
of unit spacelike ($+$) and timelike ($-$) tangent vectors; we refer to \cite{refBBZ,GKV02} for a further
discussion in the
pseudo-Riemannian context. 

In this paper, we will focus our attention on a different algebraic property of the Jacobi operator.
One says $\mathcal{M}$ is {\it Jacobi--Tsankov} if
$\mathcal{J}(x)\mathcal{J}(y)=\mathcal{J}(y)\mathcal{J}(x)$ for all tangent vectors $x,y$. One says
$\mathcal{J}$ is {\it $2$-step Jacobi nilpotent} if $\mathcal{J}(x)\mathcal{J}(y)=0$ for all tangent
vectors $x,y$. The notation is motivated by the work of Tsankov \cite{T05}.

It is convenient to work in the algebraic setting.

\begin{definition}\label{defn-1.1}
\rm Let $V$ be a finite dimensional vector space.\begin{enumerate}
\item  $A\in\otimes^4V^*$ is an {\it algebraic curvature tensor} if $A$ has the symmetries of
the curvature tensor:
\begin{enumerate}
\item $A(v_1,v_2,v_3,v_4)=-A(v_2,v_1,v_3,v_4)=A(v_3,v_4,v_1,v_2)$.
\item $A(v_1,v_2,v_3,v_4)+A(v_2,v_3,v_1,v_4)+A(v_3,v_1,v_2,v_4)=0$.
\end{enumerate}
Let $\mathfrak{A}(V)$ be the set of algebraic curvature
tensors on $V$. \item 
$\mathfrak{M}:=(V,\langle\cdot,\cdot\rangle,A)$ is a $0$-model if
$\langle\cdot,\cdot\rangle$ is a non-degenerate inner product of
signature $(p,q)$ on $V$ and if $A\in\mathfrak{A}(V)$. The associated {\it skew-symmetric curvature operator}
$\mathcal{A}(x,y)$ is characterized by
\par\centerline{$\langle\mathcal{A}(x,y)z,w\rangle=A(x,y,z,w)$,}\par\noindent
and the associated {\it Jacobi operator} is given by
$\mathcal{J}(x)y:=\mathcal{A}(y,x)x$.
\item $\mathfrak{M}$ is {\it Jacobi-Tsankov} if
$\mathcal{J}(x)\mathcal{J}(y)=\mathcal{J}(y)\mathcal{J}(x)$ $\forall$ $x,y\in V$. \item 
$\mathfrak{M}$ is {\it $2$-step Jacobi nilpotent} if $\mathcal{J}(x)\mathcal{J}(y)=0$
$\forall$ $x,y\in V$. \item  $\mathfrak{M}$ is {\it skew-Tsankov} if
$\mathcal{A}(x_1,x_2)\mathcal{A}(x_3,x_4)=\mathcal{A}(x_3,x_4)\mathcal{A}(x_1,x_2)$ $\forall$ $x_i\in V$.
\item $\mathfrak{M}$ is {\it
$2$-step skew-curvature nilpotent} if
$\mathcal{A}(x_1,x_2)\mathcal{A}(x_3,x_4)=0$ $\forall$ $x_i\in V$.
\item $\mathfrak{M}$ is {\it mixed-Tsankov} if
$\mathcal{A}(x_1,x_2)\mathcal{J}(x_3)=\mathcal{J}(x_3)\mathcal{A}(x_1,x_2)$ $\forall$ $x_i\in V$.
\item$\mathfrak{M}$ is {\it mixed-nilpotent-Tsankov} if
$\mathcal{A}(x_1,x_2)\mathcal{J}(x_3)=\mathcal{J}(x_3)\mathcal{A}(x_1,x_2)=0$ $\forall$ $x_i\in V$.
\item The {\it $0$-model of $\mathcal{M}$ at $P\in M$} is given by
$\mathfrak{M}(\mathcal{M},P):=(T_PM,g_P,R_P)$.
\item We say that $\mathcal{M}$ is a {\it geometric realization} of $\mathfrak{M}$ and that
$\mathcal{M}$ is {\it $0$-curvature homogeneous with model $\mathfrak{M}$} if for any point $P\in M$, 
$\mathfrak{M}(\mathcal{M},P)$ is isomorphic to $\mathfrak{M}$, i.e. if there exists an isomorphism
$\Theta_P:T_PM\rightarrow V$ so that $\Theta_P^*\{\langle\cdot,\cdot\rangle\}=g_P$ and so that
$\Theta_P^*A=R_P$.
\end{enumerate}\end{definition}

The following results relate these
concepts in the algebraic setting. They show in particular that any
Jacobi--Tsankov Riemannian $(p=0)$ or Lorentzian $(p=1)$ manifold is flat:

\begin{theorem}\label{thm-1.2}
Let $\mathfrak{M}:=(V,\langle\cdot,\cdot\rangle,A)$ be a $0$-model.
\begin{enumerate}
\item Let $\mathfrak{M}$ be either Jacobi--Tsankov or mixed-Tsankov. Then one has that $\mathcal{J}(x)^2=0$. Furthermore, if
$p=0$ or if $p=1$, then $A=0$.
\item If $\mathfrak{M}$ is Jacobi--Tsankov and if $\dim(V)<14$, $\mathfrak{M}$ is $2$-step Jacobi
nilpotent.
\item The following conditions are equivalent if $\mathfrak{M}$ is indecomposable:
\begin{enumerate}
\item $\mathfrak{M}$ is $2$-step Jacobi nilpotent. \item
$\mathfrak{M}$ is $2$-step skew-curvature nilpotent. \item
We can decompose $V=W\oplus\bar W$ for $W$ and $\bar W$  totally
isotropic subspaces of $V$ and for $A=A_W\oplus 0$ where the tensor
$A_W\in\mathfrak{A}(W)$ is indecomposable.
\end{enumerate}\end{enumerate}
\end{theorem}

Theorem \ref{thm-1.2} is sharp. There is a $14$-dimensional model $\mathfrak{M}_{14}$ which
is Jacobi--Tsankov but which is not $2$-step Jacobi nilpotent. This example will form the focus of
our investigations in this paper. It may be defined as follows; it is essentially unique up to
isomorphism.

\begin{definition}\label{defn-1.3}
\rm Let $\{\alpha_i,\alpha_i^*,\beta_{i,1},\beta_{i,2},\beta_{4,1},\beta_{4,2}\}$ be a basis for
$\mathbb{R}^{14}$ for
$1\le i\le 3$. Let
$\mathfrak{M}_{14}:=(\mathbb{R}^{14},\langle\cdot,\cdot\rangle,A)$ be the $0$-model where the non-zero
components of $\langle\cdot,\cdot\rangle$ and of $A$ are given, up to the usual symmetries, by:
\begin{equation}\label{eqn-1.a}
\begin{array}{l}
\langle\alpha_i,\alpha_i^*\rangle=\langle\beta_{i,1},\beta_{i,2}\rangle=1\text{ for }1\le i\le 3,\\
\langle\beta_{4,1},\beta_{4,1}\rangle=\langle\beta_{4,2},\beta_{4,2}\rangle=-\textstyle\frac12,\quad
  \langle\beta_{4,1},\beta_{4,2}\rangle=\textstyle\frac14,\gronk\\
A(\alpha_2,\alpha_1,\alpha_1,\beta_{2,1})=A(\alpha_3,\alpha_1,\alpha_1,\beta_{3,1})=1,\gronk\\
A(\alpha_3,\alpha_2,\alpha_2,\beta_{3,2})=A(\alpha_1,\alpha_2,\alpha_2,\beta_{1,2})=1,\gronk\\
A(\alpha_1,\alpha_3,\alpha_3,\beta_{1,1})=A(\alpha_2,\alpha_3,\alpha_3,\beta_{2,2})=1,\gronk\\
A(\alpha_1,\alpha_2,\alpha_3,\beta_{4,1})=A(\alpha_1,\alpha_3,\alpha_2,\beta_{4,1})=-\textstyle\frac12,\gronk\\
A(\alpha_2,\alpha_3,\alpha_1,\beta_{4,2})=A(\alpha_2,\alpha_1,\alpha_3,\beta_{4,2})=-\textstyle\frac{1}{2}\,.\gronk
\end{array}\end{equation}\end{definition}
Let $\mathcal{G}(\mathfrak{M}_{14})$ be the group of symmetries of the model:
$$
\mathcal{G}(\mathfrak{M}_{14})=\{T\in GL(\mathbb{R}^{14}):
T^*\{\langle\cdot,\cdot\rangle\}=\langle\cdot,\cdot\rangle,\ T^*A=A\}\,.
$$
Let
$SL_\pm(3):=\{A\in GL(\mathbb{R}^3):\det(A)=\pm1\}$.
In Section \ref{sect-2}, we will establish:

\begin{theorem}\label{thm-1.4}
Let $\mathfrak{M}_{14}$ be as in Definition \ref{defn-1.3}.
\begin{enumerate}\item $\mathfrak{M}_{14}$ is Jacobi--Tsankov of signature $(8,6)$.
\item $\mathfrak{M}_{14}$ is neither $2$-step Jacobi nilpotent nor skew-Tsankov.
\item There is a short exact sequence
$1\rightarrow
\mathbb{R}^{21}\rightarrow\mathcal{G}(\mathfrak{M}_{14})
\rightarrow SL_\pm(3)\rightarrow1$.
\item $\mathfrak{M}_{14}$ is mixed--Tsankov.
\end{enumerate}
\end{theorem}

In Section \ref{sect-3}, we will show that the model $\mathfrak{M}_{14}$ is geometrically
realizable. Thus there exist Jacobi--Tsankov manifolds which are not $2$-step Jacobi nilpotent.
We introduce the following notation.
\begin{definition}\label{defn-1.5}
\rm Let
$\{x_i,x_i^*,y_{i,1},y_{i,2},y_{4,1},y_{4,2}\}$ for $1\le i\le 3$ be coordinates
on $\mathbb{R}^{14}$. Suppose given a collection of functions $\Phi:=\{\phi_{i,j}\}\in
C^\infty(\mathbb{R})$ such that $\phi_{i,1}^\prime\phi_{i,2}^\prime=1$. Let
$\mathcal{M}_\Phi:=(\mathbb{R}^{14},g_\Phi)$ where the non-zero
components of $g_\Phi$ are, up to the usual $\mathbb{Z}_2$
symmetry, given by:
$$\begin{array}{ll}
g_\Phi(\partial_{x_i},\partial_{x_i^*})=g_\Phi(\partial_{y_{i,1}},\partial_{y_{i,2}})=1
   \text{ for }1\le i\le 3,\\
g_\Phi(\partial_{y_{4,1}},\partial_{y_{4,1}})=g_\Phi(\partial_{y_{4,2}},\partial_{y_{4,2}})
   =-\textstyle\frac12,&
  g_\Phi(\partial_{y_{4,1}},\partial_{y_{4,2}})=\textstyle\frac14,\vphantom{\vrule height 11pt}\\
  g_\Phi(\partial_{x_1},\partial_{x_1})=-2\phi_{2,1}(x_2)y_{2,1}-2\phi_{3,1}(x_3)y_{3,1},&
  g_\Phi(\partial_{x_2},\partial_{x_3})=x_1y_{4,1},\vphantom{\vrule height 11pt}\\
  g_\Phi(\partial_{x_2},\partial_{x_2})=-2\phi_{3,2}(x_3)y_{3,2}-2\phi_{1,2}(x_1)y_{1,2},&
g_\Phi(\partial_{x_1},\partial_{x_3})=x_2y_{4,2},\vphantom{\vrule height 11pt}\\
  g_\Phi(\partial_{x_3},\partial_{x_3})=-2\phi_{1,1}(x_1)y_{1,1}-2\phi_{2,2}(x_2)y_{2,2}
  \vphantom{\vrule height 11pt}\,.
\end{array}$$
\end{definition}

\begin{theorem}\label{thm-1.6}
Let $\mathcal{M}_\Phi:=(\mathbb{R}^{14},g_\Phi)$ be as in Definition \ref{defn-1.5}.
\begin{enumerate}
\item $\mathcal{M}_\Phi$ is geodesically complete. 
\item For all $P\in\mathbb{R}^{14}$, $\exp_P$ is a diffeomorphism from $T_P(\mathbb{R}^{14})$
to $\mathbb{R}^{14}$. 
\item  $\mathcal{M}_\Phi$ has $0$-model $\mathfrak{M}_{14}$.
\item $\mathcal{M}_\Phi$ is Jacobi--Tsankov but $\mathcal{M}_\Phi$ is not $2$-step Jacobi nilpotent.
\end{enumerate}
\end{theorem}

If we specialize the construction, we can say a bit more. We will establish the following result
in Section \ref{sect-4}:

\begin{theorem}\label{thm-1.7}
Set $\phi_{2,1}(x_2)=\phi_{2,2}(x_2)=x_2$ and $\phi_{3,1}(x_3)=\phi_{3,2}(x_3)=x_3$ in Definition
\ref{defn-1.5}. Let $\{\phi_{1,1},\phi_{1,2}\}$ be real analytic with $\phi_{1,1}^\prime\phi_{1,2}^\prime=1$ and with
$\phi_{1,j}^{\prime\prime}\ne0$. Then
\begin{enumerate}
\item
$\Xi:=\{1-\phi_{1,1}^\prime\phi_{1,1}^{\prime\prime\prime}(\phi_{1,1}^{\prime\prime})^{-2}\}^2$
is a local isometry invariant of $\mathcal{M}_\Phi$.
\item If $\phi_{1,1}^\prime(x_1)\ne be^{cx_1}$, then $\Xi$ is not locally constant and hence
$\mathcal{M}_\Phi$ is not locally homogeneous.
\end{enumerate}
\end{theorem}

There are symmetric spaces which have model $\mathfrak{M}_{14}$.

\begin{definition}\label{defn-1.8}
\rm  Let
$\{x_i,x_i^*,y_{i,1},y_{i,2},y_{4,1},y_{4,2}\}$ for $1\le i\le 3$
be coordinates
on $\mathbb{R}^{14}$. Let $A:=\{a_{i,j}\}$ be a collection of real constants. Let
$\mathcal{M}_A:=(\mathbb{R}^{14},g_A)$ where the
non-zero components of $g_A$ are given, up to the usual
$\mathbb{Z}_2$ symmetry, by:
\begin{eqnarray*}
&&{g_A}(\partial_{x_i},\partial_{x_i^*})=g_A(\partial_{y_{i,1}},\partial_{y_{i,2}})=1\text{ for }1\le
i\le 3,\\
&&{g_A}(\partial_{y_{4,1}},\partial_{y_{4,1}})=g_A(\partial_{y_{4,2}},\partial_{y_{4,2}})=-\textstyle\frac12,\qquad
  {g_A}(\partial_{y_{4,1}},\partial_{y_{4,2}})=\textstyle\frac14,\\
&&{g_A}(\partial_{x_1},\partial_{x_1})=-2 a_{2,1} x_2 y_{2,1}-2 a_{3,1}x_3 y_{3,1},\\
&&{g_A}(\partial_{x_2},\partial_{x_2})=-2 a_{3,2}x_3 y_{3,2}-2 a_{1,2}x_1 y_{1,2},\\
&&{g_A}(\partial_{x_3},\partial_{x_3})=-2 a_{1,1}x_1 y_{1,1}-2 a_{2,2}x_2 y_{2,2},\\
&&{g_A}(\partial_{x_1},\partial_{x_2})=2 (1-a_{2,1}) x_1
y_{2,1}+2(1-a_{1,2}) x_2 y_{1,2}\\
&&{g_A}(\partial_{x_2},\partial_{x_3})=x_1y_{4,1}+2 (1-a_{3,2})x_2 y_{3,2}+2(1-a_{2,2})x_3 y_{2,2},\\
&&{g_A}(\partial_{x_1},\partial_{x_3})=x_2y_{4,2}+2(1-a_{3,1})x_1
y_{3,1}+2(1-a_{1,1})x_3y_{1,1}\,.
\end{eqnarray*}
\end{definition}

We will establish the following result in Section \ref{sect-5}:

\begin{theorem}\label{thm-1.9}
Let $\mathcal{M}_A$ be described by Definition
\ref{defn-1.8}. Then $\mathcal{M}_A$ has $0$-model $\mathfrak{M}_{14}$. 
Furthermore $\mathcal{M}_A$ is
locally symmetric if
and only if\begin{enumerate}\item $a_{1,1}+a_{2,2}+a_{3,1}a_{3,2}=2$.\item
   $3 a_{2,1}+3a_{3,1}+3 a_{1,2} a_{1,1}=4$.\item
$3 a_{1,2}+3a_{3,2}+3 a_{2,1} a_{2,2}=4$.\end{enumerate}
\end{theorem}

\section{The model $\mathfrak{M}_{14}$}\label{sect-2}

We study the algebraic properties of the model $\mathfrak{M}_{14}$. Introduce the
polarization
$$\mathcal{J}(x_1,x_2):y\rightarrow\textstyle\frac12(\mathcal{A}(y,x_1)x_2+\mathcal{A}(y,x_2)x_1)\,.$$Let $\{\beta_\nu\}$ be
an enumeration of $\{\beta_{i,j}\}_{1\le i\le 4,1\le j\le 2}$. The following spaces are
invariantly defined:
\begin{eqnarray*}
&&V_{\beta,\alpha^*}:=\operatorname{Span}_{\xi_i\in
\mathbb{R}^{14}}\{\mathcal{J}(\xi_1)\xi_2\}=\operatorname{Span}\{\beta_\nu,\alpha_i^*\},\\
&&V_{\alpha^*}:=\operatorname{Span}_{\xi_i\in
\mathbb{R}^{14}}\{\mathcal{J}(\xi_1)\mathcal{J}(\xi_2)\xi_3\}=\operatorname{Span}\{\alpha_i^*\}\,.
\end{eqnarray*}

\medbreak\noindent{\it Proof of Theorem \ref{thm-1.2}}.
We have $\mathcal{J}(x)=\mathcal{J}(x,x)$ and $\mathcal{J}(x,y)x=-\frac12\mathcal{J}(x)y$. If
$\mathfrak{M}$ is Jacobi--Tsankov, then
$\mathcal{J}(x_1,x_2)\mathcal{J}(x_3,x_4)=\mathcal{J}(x_3,x_4)\mathcal{J}(x_1,x_2)$ for all $x_i$. We
may show $\mathcal{J}(x)^2=0$ by computing:
$$
0=\mathcal{J}(x,y)\mathcal{J}(x)x=\mathcal{J}(x)\mathcal{J}(x,y)x=-\textstyle\frac12\mathcal{J}(x)\mathcal{J}(x)y\,.
$$
Similarly, suppose that $\mathfrak{M}$ is mixed--Tsankov, i.e.
$$\mathcal{A}(x_1,x_2)\mathcal{J}(x_3)=\mathcal{J}(x_3)\mathcal{A}(x_1,x_2)$$ for all $x_i\in V$. We show
$\mathcal{J}(x)^2=0$ in this setting as well by computing:
$$
0=\mathcal{A}(x,y)\mathcal{J}(x)x=\mathcal{J}(x)\mathcal{A}(x,y)x=-\mathcal{J}(x)\mathcal{J}(x)y\,.
$$

We  have shown that if $\MM$ is either Jacobi-Tsankov or mixed-Tsankov, then $J(x)^2=0$. Since the Jacobi operator is
nilpotent, $\{0\}$ is the only eigenvalue of $\mathcal{J}$ so
$\MM$ is Osserman. If $p=0$, then ${\mathcal{J}}(x)$ is diagonalizable. Thus
${\mathcal{J}}(x)^2=0$ implies ${\mathcal{J}}(x)=0$ for all $x$ so $A=0$. If $p=1$, then $\MM$ is Osserman so $\MM$ has
constant sectional curvature \cite{BBG97,GKV97}; ${\mathcal{J}}(x)^2=0$, $A=0$. This establishes Assertion (1).
Assertions (2) and (3) of Theorem \ref{thm-1.2} follow from results in \cite{BG05a}.\hfill\qedbox

\medbreak \noindent{\it Proof of Theorem \ref{thm-1.4} (1,2).}
It is immediate from the definition that 
$$\mathcal{J}(\alpha_3)\mathcal{J}(\alpha_2)\alpha_1=
  \mathcal{J}(\alpha_3)\beta_{1,1}=\alpha_1^*$$
so $\mathfrak{M}_{14}$ is not $2$-step Jacobi nilpotent. 

We define $\beta_{4,1}^*$ and $\beta_{4,2}^*$ by the relations:
$\langle\beta_{4,i}^*,\beta_{4,j}\rangle=\delta_{ij}$.
We then have:
$$
\textstyle\beta_{4,1}^*=-\frac{8}{3}\beta_{4,1}-\frac{4}{3}\beta_{4,2},\quad
\textstyle\beta_{4,2}^*=-\frac{4}{3}\beta_{4,1}-\frac{8}{3}\beta_{4,2}\,.
$$
Let $\mathcal{A}_{ij}:=\mathcal{A}(\alpha_i,\alpha_j)$. We show that $\mathfrak{M}_{14}$ is not skew-Tsankov by
computing:
\begin{eqnarray*}
&&\mathcal{A}_{12}\mathcal{A}_{13}\alpha_3
 =\mathcal{A}_{12}\beta_{1,2}=-\alpha_2^*,\\
&&\mathcal{A}_{13}\mathcal{A}_{12}\alpha_3
=-\textstyle\frac12\mathcal{A}_{13}\{\beta_{4,1}^*-\beta_{4,2}^*\}
=\mathcal{A}_{13}\{\textstyle\frac23\beta_{4,1}-\textstyle\frac23\beta_{4,2}\}
  =\frac13\alpha_2^*\,.
\end{eqnarray*}

If $\xi\in \mathbb{R}^{14}$, then 
$\mathcal{J}(\xi)\alpha_i\subset V_{\beta,\alpha^*}$,
$\mathcal{J}(\xi)\beta_\nu\subset V_{\alpha^*}$, and
$\mathcal{J}(\xi)\alpha_i^*=0$.
Thus to show $\mathcal{J}(x)\mathcal{J}(y)=\mathcal{J}(y)\mathcal{J}(x)$ for all $x,y$, it suffices to show
$$\mathcal{J}(x)\mathcal{J}(y)\alpha_i=\mathcal{J}(y)\mathcal{J}(x)\alpha_i$$
for all $x,y,i$. Since
$\mathcal{J}(x)\mathcal{J}(y)\alpha_i\in V_{\alpha^*}$, this can be done by establishing:
$$
\langle \mathcal{J}(x)\alpha_i,\mathcal{J}(y)\alpha_j\rangle=\langle \mathcal{J}(y)\alpha_i,\mathcal{J}(x)\alpha_j\rangle
$$
for all $x,y,i,j$. 
Since $\mathcal{J}(x_1,x_2)\alpha_i\in V_{\alpha^*}$ if
either $x_1$ or $x_2\in V_{\beta,\alpha^*}$, we may take
$x_1=\alpha_i$ and $x_2=\alpha_j$. Let $\mathcal{J}_{ijk}:=\mathcal{J}(\alpha_i,\alpha_j)\alpha_k$.
We must show:
$$\langle\mathcal{J}_{i_1i_2i_3},\mathcal{J}_{j_1j_2j_3}\rangle=
  \langle\mathcal{J}_{i_1i_2j_3},\mathcal{J}_{j_1j_2i_3}\rangle\quad
\forall i_1i_2i_3j_1j_2j_3\,.$$

The non-zero components of $\mathcal{J}_{ijk}=\mathcal{J}_{jik}$ are:
$$\begin{array}{l}
\begin{array}{llll}
\mathcal{J}_{112}=\beta_{2,2},&\mathcal{J}_{113}=\beta_{3,2},&
\mathcal{J}_{221}=\beta_{1,1},&\mathcal{J}_{223}=\beta_{3,1},\gronk\\
\mathcal{J}_{331}=\beta_{1,2},&\mathcal{J}_{332}=\beta_{2,1},&
\mathcal{J}_{121}=-\textstyle\frac12\beta_{2,2},&\mathcal{J}_{122}=-\frac12\beta_{1,1},
        \gronk\\
\mathcal{J}_{131}=-\textstyle\frac12\beta_{3,2},
        &\mathcal{J}_{133}=-\textstyle\frac12\beta_{1,2},&
\mathcal{J}_{232}=-\textstyle\frac12\beta_{3,1},
        &\mathcal{J}_{233}=-\textstyle\frac12\beta_{2,1},\gronk
\end{array}\\
\ \mathcal{J}_{132}=\textstyle\frac14\beta_{4,1}^*-\textstyle\frac12\beta_{4,2}^*=\beta_{4,2},\qquad\quad
\mathcal{J}_{231}=-\textstyle\frac12\beta_{4,1}^*+\textstyle\frac14\beta_{4,2}^*=\beta_{4,1},\gronk\\
\ \mathcal{J}_{123}=\textstyle\frac14\beta_{4,1}^*+\textstyle\frac14\beta_{4,2}^*=-\beta_{4,1}-\beta_{4,2}
   \,.\gronk
\end{array}$$
The non-zero inner products are:
$$\begin{array}{llll}
\langle \mathcal{J}_{112},\mathcal{J}_{332}\rangle=1,&\langle
\mathcal{J}_{112},\mathcal{J}_{233}\rangle=-\textstyle\frac12,&
\langle \mathcal{J}_{121},\mathcal{J}_{332}\rangle=-\textstyle\frac12,&\langle
\mathcal{J}_{121},\mathcal{J}_{233}\rangle=\textstyle\frac14,\gronk\\
\langle \mathcal{J}_{113},\mathcal{J}_{223}\rangle=1,&\langle
\mathcal{J}_{113},\mathcal{J}_{232}\rangle=-\textstyle\frac12,&
\langle \mathcal{J}_{131},\mathcal{J}_{223}\rangle=-\textstyle\frac12,&\langle
\mathcal{J}_{232},\mathcal{J}_{131}\rangle=\textstyle\frac14,\gronk\\
\langle \mathcal{J}_{221},\mathcal{J}_{331}\rangle=1,&\langle
\mathcal{J}_{221},\mathcal{J}_{133}\rangle=-\textstyle\frac12,&
\langle \mathcal{J}_{122},\mathcal{J}_{331}\rangle=-\textstyle\frac12,&\langle
\mathcal{J}_{122},\mathcal{J}_{133}\rangle=\textstyle\frac14,\gronk\\
\langle \mathcal{J}_{123},\mathcal{J}_{123}\rangle=\star,&\langle
\mathcal{J}_{123},\mathcal{J}_{132}\rangle=\textstyle\frac14,
&\langle \mathcal{J}_{123},\mathcal{J}_{231}\rangle=\textstyle\frac14,&
\langle \mathcal{J}_{132},\mathcal{J}_{132}\rangle=\star,\gronk\\
\langle\mathcal{J}_{132},\mathcal{J}_{231}\rangle=\textstyle\frac14,&
\langle \mathcal{J}_{231},\mathcal{J}_{231}\rangle=\star
\,.\gronk
\end{array}
$$
The desired symmetries are now immediate:
\begin{eqnarray*}
&&\langle\mathcal{J}_{112},\mathcal{J}_{233}\rangle=-\textstyle\frac12=
\langle\mathcal{J}_{113},\mathcal{J}_{232}\rangle,\quad
\langle\mathcal{J}_{123},\mathcal{J}_{132}\rangle=
\textstyle\frac14=\langle \mathcal{J}_{122},\mathcal{J}_{133}\rangle,\\
&&\langle\mathcal{J}_{121},\mathcal{J}_{332}\rangle=-\textstyle\frac12=
\langle\mathcal{J}_{122},\mathcal{J}_{331}\rangle,\quad
\langle \mathcal{J}_{123},\mathcal{J}_{231}\rangle
=\textstyle\frac14=\langle
\mathcal{J}_{121},\mathcal{J}_{233}\rangle,\\
&&\langle\mathcal{J}_{131},\mathcal{J}_{223}\rangle=-\textstyle\frac12=\langle
\mathcal{J}_{133},\mathcal{J}_{221}\rangle,\quad\langle
\mathcal{J}_{132},\mathcal{J}_{231}\rangle= \textstyle\frac14
=\langle \mathcal{J}_{131},\mathcal{J}_{232}\rangle\,.\quad\qedbox
\end{eqnarray*}

\bigbreak\noindent{\it Proof of Theorem \ref{thm-1.4} (3,4).} Let $\mathcal{G}=\mathcal{G}(\mathfrak{M}_{14})$ be the
group of symmetries of the model $\mathfrak{M}_{14}$. Note that the spaces $V_{\beta,\alpha^*}$ and $V_{\alpha^*}$ are
 preserved by $\mathcal{G}$, i.e. 
\begin{equation}\label{eqn-2.a}
T V_{\alpha^*}\subset V_{\alpha^*}\quad\text{and}\quad
TV_{\beta,\alpha^*}\subset V_{\beta,\alpha^*}\quad\text{if}\quad T\in\mathcal{G}\,.
\end{equation}
Let $\tau:\mathcal{G}\rightarrow GL(3)$ be the restriction of $T$ to $V_{\alpha^*}=\mathbb{R}^3$. We will
prove Theorem \ref{thm-1.4} (3) by showing:
$$
SL_\pm(3)=\tau(\mathcal{G})\quad\text{and}\quad\ker(\tau)=\mathbb{R}^{21}\,.
$$

We argue as follows to show $SL_\pm(3)\subset\tau(\mathcal{G})$.
 Let $\beta_{4,3}:=-\beta_{4,1}-\beta_{4,2}$. One may interchange
the first two coordinates by setting:
$$\begin{array}{llll}
T:\alpha_1\leftrightarrow\alpha_2,&T:\alpha_3\leftrightarrow\alpha_3,&
T:\alpha_1^*\leftrightarrow\alpha_2^*,&
T:\alpha_3^*\leftrightarrow\alpha_3^*,\\
T:\beta_{1,1}\leftrightarrow\beta_{2,2},&
T:\beta_{1,2}\leftrightarrow\beta_{2,1},&
T:\beta_{3,1}\leftrightarrow\beta_{3,2},&
T:\beta_{4,1}\leftrightarrow\beta_{4,2}\,.
\end{array}$$
One may interchange the first and third coordinates by setting:
$$\begin{array}{llll}
T:\alpha_1\leftrightarrow\alpha_3,&T:\alpha_2\leftrightarrow\alpha_2,&T:\alpha_1^*\leftrightarrow\alpha_3^*,&
T:\alpha_2^*\leftrightarrow\alpha_2^*,\\
T:\beta_{1,1}\leftrightarrow\beta_{3,1},&
  T:\beta_{1,2}\leftrightarrow\beta_{3,2},&
T:\beta_{2,1}\leftrightarrow\beta_{2,2},&
  T:\beta_{4,1}\leftrightarrow\beta_{4,3},\\T:\beta_{4,2}\leftrightarrow\beta_{4,2}\,.
\end{array}$$
To form a rotation in the first two coordinates, we set
\vglue -.5truecm$$\begin{array}{l}
    \begin{array}{ll}
  T_\theta:\alpha_1\rightarrow \cos\theta\alpha_1+\sin\theta\alpha_2,&
  T_\theta:\alpha_2\rightarrow-\sin\theta\alpha_1+\cos\theta\alpha_2,\gronk\\
  T_\theta:\alpha_1^*\rightarrow \cos\theta\alpha_1^*+\sin\theta\alpha_2^*,&
  T_\theta:\alpha_2^*\rightarrow-\sin\theta\alpha_1^*+\cos\theta\alpha_2^*,\gronk\\
  T_\theta:\alpha_3\rightarrow\alpha_3,& T_\theta:\alpha_3^*\rightarrow\alpha_3^*,\gronk\\
  T_\theta:\beta_{1,1}\rightarrow\cos\theta\beta_{1,1}+\sin\theta\beta_{2,2},&
  T_\theta:\beta_{1,2}\rightarrow \cos\theta\beta_{1,2}+\sin\theta\beta_{2,1},\gronk\\
  T_\theta:\beta_{2,1}\rightarrow -\sin\theta\beta_{1,2}+\cos\theta\beta_{2,1},&
  T_\theta:\beta_{2,2}\rightarrow-\sin\theta\beta_{1,1}+\cos\theta\beta_{2,2},\gronk
\end{array}\\
T_\theta:\beta_{3,1}\rightarrow
\sin^2\theta\beta_{3,2}-2\sin\theta\cos\theta\beta_{4,3}+\cos^2\theta\beta_{3,1},\gronk\\
T_\theta:\beta_{3,2}\rightarrow \cos^2\theta\beta_{3,2}+2\cos\theta\sin\theta\beta_{4,3}+\sin^2\theta\beta_{3,1},\gronk\\
T_\theta:\beta_{4,1}\rightarrow \frac12 \sin\theta\cos\theta\beta_{3,2}-\frac12\sin\theta\cos\theta\beta_{3,1}
-\sin^2\theta\beta_{4,2}+\cos^2\theta\beta_{4,1},\gronk\\
T_\theta:\beta_{4,2}\rightarrow \frac12
\sin\theta\cos\theta\beta_{3,2}-\frac12\sin\theta\cos\theta\beta_{3,1}
+\cos^2\theta\beta_{4,2}-\sin^2\theta\beta_{4,1}\,.\gronk \end{array}$$
Finally, we show that
the dilatations of determinant
$1$ belong to
$\operatorname{Range}\{\tau\}$. Suppose $a_1a_2a_3=1$. We set
$$\begin{array}{llll}
T\alpha_1=a_1\alpha_1,& T\alpha_2=a_2\alpha_2,& T\alpha_3=a_3\alpha_3,& T\alpha_1^*=\textstyle\frac1{a_1}\alpha_1^*,\\
T\alpha_2^*=\textstyle\frac1{a_2}\alpha_2^*,& T\alpha_3^*=\textstyle\frac1{a_3}\alpha_3^*,&
T\beta_{1,1}=\textstyle\frac{a_2}{a_3}\beta_{1,1},&
T\beta_{1,2}=\textstyle\frac{a_3}{a_2}\beta_{1,2},\gronk\\
T\beta_{2,1}=\textstyle\frac{a_3}{a_1}\beta_{2,1},&
T\beta_{2,2}=\textstyle\frac{a_1}{a_3}\beta_{2,2},&
T\beta_{3,1}=\textstyle\frac{a_2}{a_1}\beta_{3,1},&
T\beta_{3,2}=\textstyle\frac{a_1}{a_2}\beta_{3,2},\gronk\\
T\beta_{4,1}=\beta_{4,1},&T\beta_{4,2}=\beta_{4,2}\,.\gronk
\end{array}$$
Since these elements acting on $V_{\alpha^*}$ generate $SL_\pm(3)$,  $SL_\pm(3)\subset\tau(\mathcal{G})$. Conversely,
let $T\in\mathcal{G}$. We must show $\tau(T)\in SL_\pm(3)$. As $SL_\pm(3)\subset\operatorname{Range}(\tau)$, there
exists
$S\in\mathcal{G}$ so that
$\tau(TS)$ is diagonal. Thus without loss of generality, we may assume
$\tau(T)$ is diagonal and hence:
$$
T\alpha_i=a_i\alpha_i+\textstyle\sum_\nu b_i^\nu\beta_\nu+\textstyle\sum_jc_i^j\alpha_j^*,\quad
T\beta_\nu=b_\nu\beta_\nu+\textstyle\sum_id_\nu^i\alpha_i^*,\quad
T\alpha_i^*=a_i^{-1}\alpha_i^*\,.
$$
The relations
\begin{eqnarray*}
&&\textstyle -\frac12=A(T\alpha_1,T\alpha_2,T\alpha_3,T\beta_{4,1})=\textstyle -\frac12
a_1a_2a_3b_{4,1},\\ &&-\textstyle\frac12=\langle
T\beta_{4,1},T\beta_{4,1}\rangle=-\textstyle\frac12b_{4,1}b_{4,1}
\end{eqnarray*}
show that $b_{4,1}^2=1$ and thus $a_1a_2a_3=\pm1$. Thus Range$(\tau)=SL_\pm(3)$.

We complete the proof of Assertion (3) by studying $\ker(\tau)$. If one has $T\in\ker(\tau)$, then
$$
T\alpha_i=\alpha_i+\textstyle\sum_\nu b_i^\nu\beta_\nu+\textstyle\sum_jc_i^j\alpha_j^*,\quad
 T\beta_\nu=\beta_\nu+\textstyle\sum_id_\nu^i\alpha_i^*,\quad T\alpha_i^*=\alpha_i^*\,.
$$
Using the relations $A(\alpha_i,\alpha_j,\alpha_k,\alpha_l)=0$
then leads to the following 6 linear equations the
coefficients $b_i^\nu$ must satisfy:
\medbreak\quad
$0=A(T\alpha_2,T\alpha_1,T\alpha_1,T\alpha_2)$\smallbreak\quad$\phantom{0}
=2A(b_2^{2,1}\beta_{2,1},\alpha_1,\alpha_1,\alpha_2)
  +2A(b_1^{1,2}\beta_{1,2},\alpha_2,\alpha_2,\alpha_1)
  =2b_2^{2,1}+2b_1^{1,2},$\goodbreak\smallbreak\quad$
0=A(T\alpha_3,T\alpha_1,T\alpha_1,T\alpha_3)$\smallbreak\quad$\phantom{0}
=2A(b_3^{3,1}\beta_{3,1},\alpha_1,\alpha_1,\alpha_3)
+2A(b_1^{1,1}\beta_{1,1},\alpha_3,\alpha_3,\alpha_1)
  =2b_3^{3,1}+2b_1^{1,1},$\goodbreak\par\quad$
0=A(T\alpha_3,T\alpha_2,T\alpha_2,T\alpha_3)$\smallbreak\quad$\phantom{0}
=2A(b_3^{3,2}\beta_{3,2},\alpha_2,\alpha_2,\alpha_3)
+2A(b_2^{2,2}\beta_{2,2},\alpha_3,\alpha_3,\alpha_2)
=2b_3^{3,2}+2b_2^{2,2},$\goodbreak\smallbreak\quad$
0=A(T\alpha_2,T\alpha_1,T\alpha_1,T\alpha_3)$\smallbreak\quad$\phantom{0}
=A(b_2^{3,1}\beta_{3,1},\alpha_1,\alpha_1,\alpha_3)
+A(\alpha_2,\alpha_1,\alpha_1,b_3^{2,1}\beta_{2,1})$\smallbreak\quad$
\quad+A(\alpha_2,b_1^{4,1}\beta_{4,1}+b_1^{4,2}\beta_{4,2},\alpha_1,\alpha_3)
     +A(\alpha_2,\alpha_1,b_1^{4,1}\beta_{4,1}+b_1^{4,2}\beta_{4,2},\alpha_3)$\smallbreak\quad$\phantom{0}
=b_2^{3,1}+b_3^{2,1}-\frac12b_1^{4,1}-\frac12b_1^{4,1}+\frac12b_1^{4,2},$\goodbreak\smallbreak\quad$
0=A(T\alpha_1,T\alpha_2,T\alpha_2,T\alpha_3)$\smallbreak\quad$\phantom{0}
=A(b_1^{3,2}\beta_{3,2},\alpha_2,\alpha_2,\alpha_3)
+A(\alpha_1,\alpha_2,\alpha_2,b_3^{1,2}\beta_{1,2})$\smallbreak\quad$
\quad+A(\alpha_1,b_2^{4,1}\beta_{4,1}+b_2^{4,2}\beta_{4,2},\alpha_2,\alpha_3)
     +A(\alpha_1,\alpha_2,b_2^{4,1}\beta_{4,1}+b_2^{4,2}\beta_{4,2},\alpha_3)$\smallbreak\quad$\phantom{0}
=b_1^{3,2}+b_3^{1,2}-\frac12b_2^{4,2}+\frac12b_2^{4,1}-\frac12b_2^{4,2},$\goodbreak\smallbreak\quad$
0=A(T\alpha_1,T\alpha_3,T\alpha_3,T\alpha_2)$\smallbreak\quad$\phantom{0}
=A(b_1^{2,2}\beta_{2,2},\alpha_3,\alpha_3,\alpha_2)
  +A(\alpha_1,\alpha_3,\alpha_3,b_2^{1,1}\beta_{1,1})$\smallbreak\quad$
\quad
  +A(\alpha_1,b_3^{4,1}\beta_{4,1}+b_3^{4,2}\beta_{4,2},\alpha_3,\alpha_2)
  +A(\alpha_1,\alpha_3,b_3^{4,1}\beta_{4,1}+b_3^{4,2}\beta_{4,2},\alpha_2)$\smallbreak\quad$\phantom{0}
=b_1^{2,2}+b_2^{1,1}+\frac12b_3^{4,2}+\frac12b_3^{4,1}$.
\medbreak\noindent
These equations are linearly independent so there are $18$ degrees of freedom in choosing the $b$'s. Once the $b$'s are
known, the coefficients
$d_\nu^i$ are determined;
$$0=\langle T\alpha_i,T\beta_\nu\rangle=d_\nu^i+\textstyle\sum_\mu\langle\beta_\nu,\beta_\mu\rangle b_i^\mu\,.$$
The relation $\langle T\alpha_i,T\alpha_j\rangle=\delta_{ij}$ implies $c_i^j+c_j^i=0$; this creates an additional $3$
degrees of freedom. Thus $\ker(\tau)$ is isomorphic to the additive
group
$\mathbb{R}^{21}$.

Let $\xi_i\in V$. Since $\mathcal{R}(\xi_1,\xi_2)\mathcal{J}(\xi_3)=\mathcal{J}(\xi_3)\mathcal{R}(\xi_1,\xi_2)=0$
if any of the $\xi_i\in V_{\beta,\alpha^*}$, we may work modulo $V_{\beta,\alpha^*}$ and suppose that
$\xi_i\in\operatorname{Span}\{\alpha_i\}$. Since $\mathcal{R}(\xi_1,\xi_2)=0$ if the $\xi_i$ are linearly
dependent, we suppose $\xi_1$ and $\xi_2$ are linearly independent.

There are 2 cases to be considered. We first suppose
$\xi_3\in\operatorname{Span}\{\xi_1,\xi_2\}$. The argument given above shows that a subgroup of $\mathcal{G}$ isomorphic
to $SL_\pm(3)$ acts $\operatorname{Span}\{\alpha_i\}$. Thus we may suppose
$\operatorname{Span}\{\xi_1,\xi_2\}=\operatorname{Span}\{\alpha_1,\alpha_2\}$ and that $\xi_3=\alpha_1$. Since
$\mathcal{A}(\xi_1,\xi_2)=c\mathcal{A}(\alpha_1,\alpha_2)$, we may also assume $\xi_1=\alpha_1$ and $\xi_2=\alpha_2$.
Let $\mathcal{A}_{ij}:=\mathcal{A}(\alpha_i,\alpha_i)$ and $\mathcal{J}_k:=\mathcal{J}(\alpha_k)$. We establish the
desired result by computing:
$$\begin{array}{ll}
\mathcal{A}_{12}\mathcal{J}_1\alpha_1=0,&
\mathcal{J}_1\mathcal{A}_{12}\alpha_1=-\mathcal{J}_1\beta_{2,2}=0,\\
\mathcal{A}_{12}\mathcal{J}_1\alpha_2=\mathcal{A}_{12}\beta_{2,2}=0,&
\mathcal{J}_1\mathcal{A}_{12}\alpha_2=\mathcal{J}_1\beta_{1,1}=0,\\
\mathcal{A}_{12}\mathcal{J}_1\alpha_3=\mathcal{A}_{12}\beta_{3,2}=0,&
\mathcal{J}_1\mathcal{A}_{12}\alpha_3=\textstyle\frac12\mathcal{J}_1(-\beta_{4,1}^*+\beta_{4,2}^*)=0\,.
\end{array}$$
On the other hand, if $\{\xi_1,\xi_2,\xi_3\}$ are linearly independent, we can apply a symmetry in $\mathcal{G}$ and
rescale to assume $\xi_i=\alpha_i$. We complete the proof of Theorem
\ref{thm-1.4} by computing:
$$\begin{array}{ll}
\mathcal{A}_{12}\mathcal{J}_3\alpha_1=\mathcal{A}_{12}\beta_{1,2}=-\alpha_2^*,&
\mathcal{J}_3\mathcal{A}_{12}\alpha_1=-\mathcal{J}_3\beta_{2,2}=-\alpha_2^*,\\
\mathcal{A}_{12}\mathcal{J}_3\alpha_2=\mathcal{A}_{12}\beta_{2,1}=\alpha_1^*,&
\mathcal{J}_3\mathcal{A}_{12}\alpha_2=\mathcal{J}_3\beta_{1,1}=\alpha_1^*,\\
\mathcal{A}_{1,2}\mathcal{J}_3\alpha_3=0,&
\mathcal{J}_3\mathcal{A}_{1,2}\alpha_3=\textstyle\frac12\mathcal{J}_3(-\beta_{4,1}^*+\beta_{4,2}^*)=0\,.\qquad\qedbox
\end{array}$$

\begin{remark}\label{rmk-2.1}\rm If $\{e_1,e_2\}$ is an oriented orthonormal basis for a non-degenerate $2$-plane $\pi$,
one may define $\mathcal{R}(\pi):=\mathcal{R}(e_1,e_2)$ and  one may define
$\mathcal{J}(\pi):=\langle e_1,e_1\rangle\mathcal{J}(e_1)+\langle e_2,e_2\rangle\mathcal{J}(e_2)$.
These operators are independent of the particular orthonormal basis chosen. Stanilov and Videv \cite{SV04} have shown
that if
$\mathcal{M}$ is a $4$-dimensional Riemannian manifold, then
$\mathcal{R}(\pi)\mathcal{J}(\pi)=\mathcal{J}(\pi)\mathcal{R}(\pi)$ for all oriented $2$-planes $\pi$ if and only if
$\mathcal{M}$ is Einstein. Assertion (4) of Theorem
\ref{thm-1.4} shows $\mathfrak{M}_{14}$ has this property. 
\end{remark}

\section{A geometric realization of $\mathfrak{M}$}\label{sect-3}

We begin the proof of Theorem \ref{thm-1.6} with a general construction:
\begin{definition}\label{defn-3.1}
\rm Let $\{x_i,x^*_i,y_\mu\}$ be coordinates on $\mathbb{R}^{2a+b}$ where $1\le
i\le a$ and $1\le\mu\le b$. We suppose given a non-degenerate symmetric matrix $C_{\mu\nu}$ and
smooth functions $\psi_{ij\mu}=\psi_{ij\mu}(\vec x)$ with $\psi_{ij\mu}=\psi_{ji\mu}$.
Consider the pseudo-Riemannian manifold $\mathcal{M}_{C,\psi}:=(\mathbb{R}^{2a+b},g_{C,\psi})$, where:
$$g_{C,\psi}(\partial_{x_i},\partial_{x_j})=2\textstyle\sum_ky_\mu\psi_{ij\mu},\quad
  g_{C,\psi}(\partial_{x_i},\partial_{x_i^*})=1,\quad
  g_{C,\psi}(\partial_{y_\mu},\partial_{y_\nu})=C_{\mu\nu}\,.$$
\end{definition}
\begin{lemma}\label{lem-3.2}
Let $\mathcal{M}_{C,\psi}=(\mathbb{R}^{2a+b},g_{C,\psi})$ be as in Definition \ref{defn-3.1}. Then
\begin{enumerate}\item $\mathcal{M}_{C,\psi}$ is geodesically complete.
\item For all $P\in\mathbb{R}^{2a+b}$, $\exp_P$ is a diffeomorphism from $T_P(\mathbb{R}^{2a+b})$
to $\mathbb{R}^{2a+b}$.
\item The possibly non-zero components of the curvature tensor are, up to the usual $\mathbb{Z}_2$
symmetries given by:
\begin{eqnarray*}
&&R(\partial_{x_i},\partial_{x_j},\partial_{x_k},\partial_{y_\nu})=-\partial_{x_i}\psi_{jk\nu}+
\partial_{x_j}\psi_{ik\nu},\\
&&R(\partial_{x_i},\partial_{x_j},\partial_{x_k},\partial_{x_l})=
\textstyle\sum_{\nu\mu}C^{\nu\mu}\left\{\psi_{ik\mu}\psi_{jl\nu}-\psi_{il\mu}\psi_{jk\nu}\right\}\\
&&\quad+\textstyle\sum_\nu y_\nu\left\{\partial_{x_i}\partial_{x_k}\psi_{jl\nu}
   +\partial_{x_j}\partial_{x_l}\psi_{ik\nu}
   -\partial_{x_i}\partial_{x_l}\psi_{jk\nu}-\partial_{x_j}\partial_{x_k}\psi_{il\nu}\right\}\,.
\end{eqnarray*}
\end{enumerate}
\end{lemma}
\begin{proof} The non-zero Christoffel symbols of the first kind are given by:
\begin{eqnarray*}
&&g(\nabla_{\partial_{x_i}}\partial_{x_j},\partial_{x_k})=
\textstyle\sum_\mu\{\partial_{x_i}\psi_{jk\mu}+\partial_{x_j}\psi_{ik\mu}-\partial_{x_k}\psi_{ij\mu}\}
y_\mu,\\ &&g(\nabla_{\partial_{x_i}}\partial_{x_j},\partial_{y_\nu})=-\psi_{ij\nu},\\
&&g(\nabla_{\partial_{x_i}}\partial_{y_\nu},\partial_{x_k})=g(\nabla_{\partial_{y_\nu}}\partial_{x_i},
\partial_{x_k})=\psi_{ik\nu},
\end{eqnarray*}
and the non-zero Christoffel symbols of the second kind are given by:
\begin{eqnarray*}
&&\nabla_{\partial_{x_i}}\partial_{x_j}=
\textstyle\sum_\mu y_\mu\{\partial_{x_i}\psi_{jk\mu}+\partial_{x_j}\psi_{ik\mu}
-\partial_{x_k}\psi_{ij\mu}\}\partial_{x_k^*}
-\textstyle\sum_{\mu\nu} C^{\nu\mu}\psi_{ij\nu}\partial_{y_{\mu}},\\
&&\nabla_{\partial_{x_i}}\partial_{y_\nu}=\nabla_{\partial_{y_\nu}}\partial_{x_i}=
\textstyle\sum_k\psi_{ik\nu}\partial_{x_k^*}\,.
\end{eqnarray*}
This shows that $\mathcal{M}$ is a generalized plane wave manifold; Assertions (1) and (2) then follow
from results in \cite{GN05}. Assertion (3) now follows by a direct calculation.
\end{proof}

\medbreak\noindent{\it Proof of Theorem \ref{thm-1.6} (1)-(3)} Assertions (1) and (2) of Theorem \ref{thm-1.6} follow by
specializing the corresponding results of Lemma \ref{lem-3.2}. We use Assertion (3) of Lemma \ref{lem-3.2} to
see that the possibly non-zero components of the curvature tensor defined by the metric of Definition \ref{defn-1.5}
are:
\medbreak\qquad
$R(\partial_{x_{i_1}},\partial_{x_{i_2}},\partial_{x_{i_3}},\partial_{x_{i_4}})=\star$,\smallbreak\qquad
$R(\partial_{x_1},\partial_{x_2},\partial_{y_{2,1}},\partial_{x_1})=\partial_{x_2}\phi_{2,1},\quad
  R(\partial_{x_1},\partial_{x_3},\partial_{y_{3,1}},\partial_{x_1})=\partial_{x_3}\phi_{3,1}$,\smallbreak\qquad
$R(\partial_{x_2},\partial_{x_3},\partial_{y_{3,2}},\partial_{x_2})=\partial_{x_3}\phi_{3,2},\quad
  R(\partial_{x_2},\partial_{x_1},\partial_{y_{1,2}},\partial_{x_2})=\partial_{x_1}\phi_{1,2}$,\smallbreak\qquad
$R(\partial_{x_3},\partial_{x_1},\partial_{y_{1,1}},\partial_{x_3})=\partial_{x_1}\phi_{1,1},\quad
  R(\partial_{x_3},\partial_{x_2},\partial_{y_{2,2}},\partial_{x_3})=\partial_{x_2}\phi_{2,2}$,\smallbreak\qquad
$R(\partial_{x_2},\partial_{x_1},\partial_{y_{4,1}},\partial_{x_3})=
  R(\partial_{x_3},\partial_{x_1},\partial_{y_{4,1}},\partial_{x_2})=-\textstyle\frac12$,\smallbreak\qquad
$R(\partial_{x_1},\partial_{x_2},\partial_{y_{4,2}},\partial_{x_3})=
  R(\partial_{x_3},\partial_{x_2},\partial_{y_{4,2}},\partial_{x_1})=-\textstyle\frac12$.\medbreak\noindent
We introduce the following basis as a first step in the proof of Assertion (3). Let
the index $i$ range from $1$ to $3$ and the index $j$ run from $1$ to $2$. Set:
\begin{equation}\label{eqn-3.a}
\bar\alpha_i:=\partial_{x_i},\quad\alpha_i^*:=\partial_{x_i^*},\quad
  \bar\beta_{4,j}:=\partial_{y_{4,j}},\quad
  \bar\beta_{i,j}:=\{\phi_{i,j}^\prime\}^{-1}\partial_{y_{i,j}}\,.
\end{equation}
Since 
$\phi_{i,1}^\prime\cdot\phi_{i,2}^\prime=1$, the
relations of Equation (\ref{eqn-1.a}) are satisfied. However, we
still have the following potentially non-zero terms to deal with:
$$g(\bar\alpha_i,\bar\alpha_j)=\star\quad\text{and}\quad
R(\bar\alpha_i,\bar\alpha_j,\bar\alpha_k,\bar\alpha_l)=\star\,.
$$
To deal with the extra curvature terms, we introduce a modified basis setting:
\begin{equation}\label{eqn-3.b}
\begin{array}{l}
\tilde\alpha_1:=\bar\alpha_1+R(\bar\alpha_1,\bar\alpha_2,\bar\alpha_3,\bar\alpha_1)\bar\beta_{4,1}
   -\textstyle\frac12R(\bar\alpha_1,\bar\alpha_2,\bar\alpha_2,\bar\alpha_1)\bar\beta_{1,2},\\
\tilde\alpha_2:=\bar\alpha_2+R(\bar\alpha_2,\bar\alpha_1,\bar\alpha_3,\bar\alpha_2)
   \bar\beta_{4,2}
-\textstyle\frac12R(\bar\alpha_2,\bar\alpha_3,\bar\alpha_3,\bar\alpha_2)\bar\beta_{2,2},\gronk\\
\tilde\alpha_3:=\bar\alpha_3-2R(\bar\alpha_3,\bar\alpha_1,\bar\alpha_2,\bar\alpha_3)
   \bar\beta_{4,1}
    -\textstyle\frac12R(\bar\alpha_1,\bar\alpha_3,\bar\alpha_3,\bar\alpha_1)\bar\beta_{3,1},\gronk\\
\beta_{1,1}:=\bar\beta_{1,1}
 +\textstyle\frac12R(\bar\alpha_1,\bar\alpha_2,\bar\alpha_2,\bar\alpha_1)\alpha_1^*,\qquad\qquad
  \beta_{1,2}:=\bar\beta_{1,2}\gronk\\
\beta_{2,1}:=\bar\beta_{2,1}
 +\textstyle\frac12R(\bar\alpha_2,\bar\alpha_3,\bar\alpha_3,\bar\alpha_2)\alpha_2^*,\qquad\qquad
  \beta_{2,2}:=\bar\beta_{2,2},\gronk\\
\beta_{3,2}:=\bar\beta_{3,2}
  +\textstyle\frac12R(\bar\alpha_1,\bar\alpha_3,\bar\alpha_3,\bar\alpha_1)\alpha_3^*,\qquad\qquad
\beta_{3,1}:=\bar\beta_{3,1},\gronk\\
\beta_{4,1}:=\bar\beta_{4,1}
  +\textstyle\frac12R(\bar\alpha_1,\bar\alpha_2,\bar\alpha_3,\bar\alpha_1)\alpha_1^*
   -\textstyle\frac14R(\bar\alpha_2,\bar\alpha_1,\bar\alpha_3,\bar\alpha_2)\alpha_2^*\gronk\\
\qquad\qquad\qquad-R(\bar\alpha_3,\bar\alpha_1,\bar\alpha_2,\bar\alpha_3)\alpha_3^*,\gronk\\
\beta_{4,2}:=\bar\beta_{4,2}
  -\textstyle\frac14R(\bar\alpha_1,\bar\alpha_2,\bar\alpha_3,\bar\alpha_1)\alpha_1^*
+\textstyle\frac12R(\bar\alpha_2,\bar\alpha_1,\bar\alpha_3,\bar\alpha_2)\alpha_2^*\gronk\\
\qquad\qquad\qquad+\textstyle\frac12R(\bar\alpha_3,\bar\alpha_1,\bar\alpha_2,\bar\alpha_3)\alpha_3^*\,.
\end{array}\end{equation}
All the normalizations of Equation (\ref{eqn-1.a}) are satisfied except for the unwanted metric terms
$g(\tilde\alpha_i,\tilde\alpha_j)$. To eliminate these terms and to exhibit a basis with the required normalizations,
we set:
\begin{equation}\label{eqn-3.c}
\alpha_i:=\tilde\alpha_i-\textstyle\frac12\sum_jg(\tilde\alpha_i,\tilde\alpha_j)\alpha_j^*\,.\qquad\qedbox
\end{equation}

\section{Isometry Invariants}\label{sect-4}

\medbreak We now turn to the task of constructing invariants. 

\begin{lemma}\label{lem-4.1}
Adopt the assumptions of Theorem \ref{thm-1.7}. Let $\{\alpha_i,\beta_\nu,\alpha_i^*\}$ be
defined by Equations   (\ref{eqn-3.a})-(\ref{eqn-3.c}). Set $\phi_1:=\phi_{1,1}^\prime$ and $\phi_2:=\phi_{1,2}^\prime$.
\begin{enumerate}
\item $\nabla R(v_1,v_2,v_3,v_4;v_5)=0$ if at least one of the $v_i\in V_{\alpha^*}$.
\item $\nabla R(v_1,v_2,v_3,v_4;v_5)=0$ if at least two of the $v_i\in V_{\beta,\alpha^*}$.
\item $\nabla^kR(\alpha_1,\alpha_2,\alpha_2,\beta_{1,2};\alpha_1,...,\alpha_1)
    =\phi_2^{-1}\phi_2^{(k)}$.
\item $\nabla^kR(\alpha_1,\alpha_3,\alpha_3,\beta_{1,1};\alpha_1,...,\alpha_1)
     =\phi_1^{-1}\phi_1^{(k)}$.
\item $\nabla R(\alpha_i,\alpha_j,\alpha_k,\beta_\nu;\alpha_{l_1},...,\alpha_{l_k})=0$ in
cases other than those given in (3) and (4) up to the usual $\mathbb{Z}_2$ symmetry
in the first 2 entries.
\end{enumerate}
\end{lemma}

\begin{proof} Let $v_i$ be coordinate vector fields. To prove Assertion (1), we suppose some $v_i\in V_{\alpha^*}$. We
may use the second Bianchi identity and the other curvature symmetries to assume without loss of generality
$v_1\in V_{\alpha^*}$. Since $\nabla_{v_5}v_1=0$ and since $R(v_1,\cdot,\cdot,\cdot)=0$, Assertion (1) follows. The proof
of the second assertion is similar and uses the fact that
$R(\cdot,\cdot,\cdot,\cdot)=0$ if $2$-entries belong to $V_{\beta,\alpha^*}$. 
The proof of the remaining assertions is similar and uses the particular form of the warping
functions $\phi_{i,j}$; the factor of $\phi_{1,j}^{-1}$ arising from the normalization in Equation
(\ref{eqn-3.a}).
\end{proof}

\begin{definition}\label{defn-4.2}
We say that a basis $\tilde{\mathcal{B}}:=\{\tilde\alpha_i,\tilde\beta_\nu,\tilde\alpha_i^*\}$ is 
{\it $0$-normalized} if the normalizations of Equation (\ref{eqn-1.a}) are satisfied and {\it
$1$-normalized} if it is $0$-normalized and if additionally
\begin{eqnarray*}
&&\nabla R(\tilde\alpha_1,\tilde\alpha_3,\tilde\alpha_3,\tilde\beta_{1,1};\tilde\alpha_1)
=-\nabla R(\tilde\alpha_3,\tilde\alpha_1,\tilde\alpha_3,\tilde\beta_{1,1};\tilde\alpha_1)
\ne0,\\
&&\nabla R(\tilde\alpha_1,\tilde\alpha_2,\tilde\alpha_2,\tilde\beta_{1,2};\tilde\alpha_1)
=-\nabla R(\tilde\alpha_2,\tilde\alpha_1,\tilde\alpha_2,\tilde\beta_{1,2};\tilde\alpha_1)
\ne0,\\
&&\nabla
R(\tilde\alpha_i,\tilde\alpha_j,\tilde\alpha_k,\tilde\beta_\nu;\tilde\alpha_l)=0\quad\text{otherwise}
\,.
\end{eqnarray*}\end{definition}

\begin{lemma}\label{lem-4.3} Adopt the assumptions of Theorem \ref{thm-1.7}. Then:
\begin{enumerate}
\item There exists a $1$-normalized basis.
\item If $\tilde{\mathcal{B}}$ is a $1$-normalized basis, then there exist constants $a_i$ so
$a_1a_2a_3=\varepsilon$ for $\varepsilon=\pm1$ and so that exactly one of the following conditions holds:
\begin{enumerate}
\item $\tilde\alpha_1=a_1\alpha_1$, $\tilde\alpha_2=a_2\alpha_2$, $\tilde\alpha_3=a_3\alpha_3$,\newline
$\tilde\beta_{1,1}=\varepsilon\textstyle\frac{a_2}{a_3}\beta_{1,1}$,
$\tilde\beta_{1,2}=\varepsilon\textstyle\frac{a_3}{a_2}\beta_{1,2}$.
\item  $\tilde\alpha_1=a_1\alpha_1$, $\tilde\alpha_2=a_3\alpha_3$, $\tilde\alpha_3=a_2\alpha_2$,\newline 
$\tilde\beta_{1,1}=\varepsilon\textstyle\frac{a_3}{a_2}\beta_{1,2}$,
$\tilde\beta_{1,2}=\varepsilon\textstyle\frac{a_2}{a_3}\beta_{1,1}$. 
\end{enumerate}
\end{enumerate}
\end{lemma}

\begin{proof} We use Equations (\ref{eqn-3.a}), (\ref{eqn-3.b}), and (\ref{eqn-3.c}) to construct a
$0$-normalized basis and then apply Lemma \ref{lem-4.1} to see this basis is $1$-normalized.
 On the other hand, if $\tilde{\mathcal{B}}$ is a $1$-normalized basis, we may
expand:
\begin{eqnarray*}
&&\tilde\alpha_1=a_{11}\alpha_1+a_{12}\alpha_2+a_{13}\alpha_3+...,\\
&&\tilde\alpha_2=a_{21}\alpha_1+a_{22}\alpha_2+a_{23}\alpha_3+...,
  \quad\tilde\beta_{1,2}=b_{21}\beta_{1,1}+b_{22}\beta_{1,2}+...\\
&&\tilde\alpha_3=a_{31}\alpha_1+a_{32}\alpha_2+a_{33}\alpha_3+...,
\quad\tilde\beta_{1,1}=b_{11}\beta_{1,1}+b_{12}\beta_{1,2}+...\,.
\end{eqnarray*}
Because 
\begin{eqnarray*}
0&\ne&\nabla R(\tilde\alpha_1,\tilde\alpha_2,\tilde\alpha_2,\tilde\beta_{1,2};\tilde\alpha_1)\\
&=&a_{11}\left\{(a_{11}a_{22}-a_{12}a_{21})a_{22}b_{22})\phi_2^{-1}\phi_2^{\prime}\right.\\
&&\quad+\left.(a_{11}a_{33}-a_{13}a_{31})a_{33}b_{21})\phi_1^{-1}\phi_1^{\prime}\right\},
\end{eqnarray*}
we have $a_{11}\ne0$. Because
\begin{eqnarray*} 
0&=&\nabla R(\tilde\alpha_1,\tilde\alpha_2,\tilde\alpha_2,\tilde\beta_{1,2};\tilde\alpha_2)
=\textstyle\frac{a_{21}}{a_{11}}\nabla
R(\tilde\alpha_1,\tilde\alpha_2,\tilde\alpha_2,\tilde\beta_{1,2};\tilde\alpha_1),
\end{eqnarray*}
we have $a_{21}=0$; similarly $a_{31}=0$. 
Since $\operatorname{Span}\{\alpha_i\}=\operatorname{Span}\{\tilde\alpha_i\}$ mod $V_{\beta,\alpha*}$,
$$a_{22}a_{33}-a_{23}a_{32}\ne0\,.$$

By hypothesis
$R(\tilde\alpha_1,\tilde\alpha_2,\tilde\alpha_3,\beta;\tilde\alpha_1)=0$ if
$\beta\in\operatorname{Span}\{\tilde\beta_\nu,\tilde\alpha_i^*\}=V_{\beta,\alpha^*}$ so
\begin{eqnarray*}
0&=&R(\tilde\alpha_1,\tilde\alpha_2,\tilde\alpha_3,\beta_{1,2};\tilde\alpha_1)=a_{11}^2a_{22}a_{32}
     \phi_2^{-1}\phi_2^{\prime},\\
0&=&R(\tilde\alpha_1,\tilde\alpha_2,\tilde\alpha_3,\beta_{1,1};\tilde\alpha_1)=a_{11}^2a_{23}a_{33}
     \phi_1^{-1}\phi_1^{\prime}\,.
\end{eqnarray*}

Suppose that $a_{22}\ne0$. Since $a_{11}^2a_{22}a_{32}=0$ and $a_{11}\ne0$, $a_{32}=0$. Since
$a_{22}a_{33}-a_{23}a_{32}\ne0$, $a_{33}\ne0$. Since $a_{11}^2a_{23}a_{33}=0$, we also have $a_{23}=0$. Since the basis
is also $0$-normalized, $diag(a_{11}^{-1},a_{22}^{-1},a_{33}^{-1})\in SL_\pm(3)$ from the discussion in Section
\ref{sect-2}. Thus
$\varepsilon:=a_{11}a_{22}a_{33}=\pm1$, $b_{11}=\varepsilon\textstyle\frac{a_{33}}{a_{22}}$, and
$b_{22}=\varepsilon\textstyle\frac{a_{22}}{a_{33}}$. These are the relations of Assertion (2a). The argument is similar
if
$a_{32}\ne0$; we simply  reverse the roles of
$\tilde\alpha_2$ and
$\tilde\alpha_3$ to establish the relations of Assertion (2b).\end{proof}

\medbreak\noindent{\it Proof of Theorem \ref{thm-1.7}.} Let
$$
\Xi(\mathcal{B}):=\frac14\left\{\frac{\nabla^2R(\tilde\alpha_1,\tilde\alpha_2,\tilde\alpha_2,\tilde\beta_{1,2};
   \tilde\alpha_1,\tilde\alpha_1)}
{\{\nabla  R(\tilde\alpha_1,\tilde\alpha_2,\tilde\alpha_2,\tilde\beta_{1,2};\tilde\alpha_1)\}^2}
-\frac{\nabla^2R(\tilde\alpha_1,\tilde\alpha_3,\tilde\alpha_3,\tilde\beta_{1,1};\tilde\alpha_1,\tilde\alpha_1)}
{\{\nabla  R(\tilde\alpha_1,\tilde\alpha_3,\tilde\alpha_3,\tilde\beta_{1,1};\tilde\alpha_1)\}^2}\right\}^2
$$

We apply Lemma \ref{lem-4.3}.
Suppose the conditions of Assertion (2a) hold. Then:
\medbreak\quad
$\nabla  R(\tilde\alpha_1,\tilde\alpha_2,\tilde\alpha_2,\tilde\beta_{1,2};\tilde\alpha_1)
    =a_1\phi_2^{-1}\phi_2^\prime$,\smallbreak\quad$
  \nabla^2R(\tilde\alpha_1,\tilde\alpha_2,\tilde\alpha_2,\tilde\beta_{1,2};\tilde\alpha_1,\tilde\alpha_1)
    =a_1^2\phi_2^{-1}\phi_2^{\prime\prime}$,\smallbreak\quad
$\nabla  R(\tilde\alpha_1,\tilde\alpha_3,\tilde\alpha_3,\tilde\beta_{1,1};\tilde\alpha_1)
    =a_1\phi_1^{-1}\phi_1^\prime$,\smallbreak\quad$ 
  \nabla^2R(\tilde\alpha_1,\tilde\alpha_3,\tilde\alpha_3,\tilde\beta_{1,1};\tilde\alpha_1,\tilde\alpha_1)
   =a_1^2\phi_1^{-1}\phi_1^{\prime\prime}$,\smallbreak\quad
$$\Xi(\mathcal{B})=\frac14\left\{\frac{\phi_2\phi_2^{\prime\prime}}{\phi_2^\prime\phi_2^\prime}-
    \frac{\phi_1\phi_1^{\prime\prime}}{\phi_1^\prime\phi_1^\prime}\right\}^2\,.
$$
The roles of $\phi_1$ and $\phi_2$ are reversed if Assertion (2b) holds. It now follows that $\Xi$ is a local isometry invariant.
Since $\phi_2=\phi_1^{-1}$, $\phi_2^\prime=-\phi_1^{-2}\phi_1^{\prime}$,
$\phi_2^{\prime\prime}=2\phi_1^{-3}\phi_1^\prime\phi_1^\prime-\phi_1^{-2}\phi_1^{\prime\prime}$,
we may establish Assertion (1) of Theorem \ref{thm-1.7} by computing
$$\frac{\phi_2\phi_2^{\prime\prime}}{\phi_2^\prime\phi_2^\prime}=
\frac{\phi_1^{-1}(2\phi_1^{-3}\phi_1^\prime\phi_1^\prime-\phi_1^{-2}\phi_1^{\prime\prime})}
{\phi_1^{-4}\phi_1^{\prime}\phi_1^{\prime}}=2-\frac{\phi_1\phi_1^{\prime\prime}}{\phi_1^\prime\phi_1^\prime}\,.
$$
Consequently
$$
\Xi=\frac14\left\{2-2\frac{\phi_1\phi_1^{\prime\prime}}{\phi_1^\prime\phi_1^\prime}\right\}^2\,.
$$

If $\mathcal{M}_\Phi$ is locally homogeneous, then $\Xi$ must be constant. Conversely, if $\Xi$ is constant, then
$\phi_1\phi_1^{\prime\prime}=k\phi_1^\prime\phi_1^\prime$ for some $k\in\mathbb{R}$.
The solutions to this ordinary differential equation take the form
$\phi_1(t)=a(t+b)^c$ if $k\ne1$ and $\phi_1(t)=ae^{bt}$ if $k=1$ for suitably chosen constants $a$ and $b$ and for
$c=c(k)$. The first family is ruled out as
$\phi_1$ and $\phi_1^\prime$ must be invertible for all $t$. Thus $\phi_1(t)$ is a pure exponential;
Assertion (2) of Theorem \ref{thm-1.7} follows.\hfill\qedbox 

\section{A symmetric space with model $\mathfrak{M}_{14}$}\label{sect-5}
We give the proof of Theorem \ref{thm-1.9} as follows. Let $\mathcal{M}_A$ be as described in Definition \ref{defn-1.8}. By
Lemma
\ref{lem-3.2} one has that:

\medbreak\quad
$R(\partial_{x_2},\partial_{x_1},\partial_{x_1},\partial_{y_{2,1}})
  =R(\partial_{x_3},\partial_{x_1},\partial_{x_1},\partial_{y_{3,1}})=1$,\smallbreak\quad
$R(\partial_{x_3},\partial_{x_2},\partial_{x_2},\partial_{y_{3,2}})
  =R(\partial_{x_1},\partial_{x_2},\partial_{x_2},\partial_{y_{1,2}})=1$,\smallbreak\quad
$R(\partial_{x_1},\partial_{x_3},\partial_{x_3},\partial_{y_{1,1}})
  =R(\partial_{x_2},\partial_{x_3},\partial_{x_3},\partial_{y_{2,2}})=1$,\smallbreak\quad
$R(\partial_{x_1},\partial_{x_2},\partial_{x_3},\partial_{y_{4,1}})
  =R(\partial_{x_1},\partial_{x_3},\partial_{x_2},\partial_{y_{4,1}})=-\textstyle\frac12$,\smallbreak\quad
$R(\partial_{x_2},\partial_{x_3},\partial_{x_1},\partial_{y_{4,2}})
  =R(\partial_{x_2},\partial_{x_1},\partial_{x_3},\partial_{y_{4,2}})=-\textstyle\frac{1}{2}$.
\medbreak The same argument constructing a $0$-normalized basis which was given in the proof  of Theorem
\ref{thm-1.4} can then be used to construct a $0$-normalized basis in this setting and establish that
$\mathcal{M}_{\mathcal{A}}$ has
$0$-model
$\mathfrak{M}_{14}$. 

We can also apply Lemma \ref{lem-3.2} to see:\medbreak\quad
$R(\partial_{x_1},\partial_{x_2},\partial_{x_2},\partial_{x_1})=-a_{3,1}a_{3,2}x_3^2$,\smallbreak\quad
$R(\partial_{x_1},\partial_{x_3},\partial_{x_3},\partial_{x_1})=-\frac{1}{3}(2+3a_{2,1}a_{2,2})x_2^2$,\smallbreak\quad
$R(\partial_{x_3},\partial_{x_2},\partial_{x_2},\partial_{x_3})=-\frac{1}{3}(2+3a_{1,1}a_{1,2})x_1^2$,\smallbreak\quad
$R(\partial_{x_2},\partial_{x_1},\partial_{x_1},\partial_{x_3})=(1-a_{1,1}-a_{1,2}+a_{1,1}a_{1,2}+a_{2,1}$\par\qquad
$-a_{2,1}a_{2,2}+a_{3,1}-a_{3,1}a_{3,2})x_2x_3$,\smallbreak\quad
$R(\partial_{x_1},\partial_{x_2},\partial_{x_2},\partial_{x_3})=(1+a_{1,2}-a_{2,1}-a_{1,1}a_{1,2}-a_{2,2}$\par\qquad
$+a_{2,1}a_{2,2}+a_{3,2}-a_{3,1}a_{3,2})x_1x_3$,\smallbreak\quad
$R(\partial_{x_1},\partial_{x_3},\partial_{x_3},\partial_{x_2})=(\frac{2}{3}+a_{1,1}-a_{1,1}a_{1,2}+a_{2,2}$\par\qquad
$-a_{2,1}a_{2,2}-a_{3,1}-a_{3,2}+a_{3,1}a_{3,2})x_1x_2$.

\medbreak 
The Christoffel symbols describing $\nabla_{\partial_{x_i}}\partial_{x_j}$ are given by:\medbreak\quad
$\nabla_{\partial_{x_1}}\partial_{x_1}=
(2-a_{2,1})y_{2,1}\partial_{x_2^*}+(2-a_{3,1})y_{3,1}\partial_{x_3^*}
    +a_{2,1}x_2\partial_{y_{2,2}}$\par\qquad$+a_{3,1}x_3\partial_{y_{3,2}}$,
\smallbreak\quad
$\nabla_{\partial_{x_2}}\partial_{x_2}=
(2-a_{1,2})y_{1,2}\partial_{x_1^*}+(2-a_{3,2})y_{3,2}\partial_{x_3^*}
+a_{1,2}x_1\partial_{y_{1,1}}$\par\qquad$+a_{3,2}x_3\partial_{y_{3,1}}$,\smallbreak\quad
$\nabla_{\partial_{x_3}}\partial_{x_3}=
(2-a_{1,1})y_{1,1}\partial_{x_1^*}+(2-a_{2,2})y_{2,2}\partial_{x_2^*}
   +a_{2,2}x_2\partial_{y_{2,1}}$\par\qquad$+a_{1,1}x_1\partial_{y_{1,2}}$,\smallbreak\quad
$\nabla_{\partial_{x_1}}\partial_{x_2}=-a_{2,1}y_{2,1}\partial_{x_1^*}-a_{1,2}y_{1,2}\partial_{x_2^*}
+\frac{y_{4,1}+y_{4,2}}{2}\partial_{x_3^*}$\par
$\qquad+(a_{1,2}-1)x_2\partial_{y_{1,1}}+(a_{2,1}-1)x_1\partial_{y_{2,2}},$\smallbreak\quad
$\nabla_{\partial_{x_1}}\partial_{x_3}=-a_{3,1}y_{3,1}\partial_{x_1^*}
+\frac{y_{4,1}-y_{4,2}}{2}\partial_{x_2^*}-a_{1,1}y_{1,1}\partial_{x_3^*}$\par\qquad
$+(a_{1,1}-1)x_3\partial_{y_{1,2}}+(a_{3,1}-1)x_1\partial_{y_{3,2}}
   +\frac{2x_2}{3}\partial_{y_{4,1}}+\frac{4x_2}{3}\partial_{y_{4,2}}$,\smallbreak\quad
$\nabla_{\partial_{x_2}}\partial_{x_3}
=\frac{-y_{4,1}+y_{4,2}}{2}\partial_{x_1^*}-a_{3,2}y_{3,2}\partial_{x_2^*}-a_{2,2}y_{2,2}\partial_{x_3^*}$
   \par\qquad
$+(a_{2,2}-1)x_3\partial_{y_{2,1}}+(a_{3,2}-1)x_2\partial_{y_{3,1}}
+\frac{4x_1}{3}\partial_{y_{4,1}}+\frac{2x_1}{3}\partial_{y_{4,2}}$.
\medbreak It is now easy to show that the non-zero components of $\nabla R$ are:
\medbreak\quad
$\nabla
R(\partial_{x_1},\partial_{x_2},\partial_{x_2},\partial_{x_1};\partial_{x_3})
   =-2(-2+a_{1,1}+a_{2,2}+a_{3,1}a_{3,2})x_3$,\smallbreak\quad
$\nabla
R(\partial_{x_1},\partial_{x_3},\partial_{x_3},\partial_{x_1};\partial_{x_2})=
-\frac{2}{3}(-4+3a_{1,2}+3a_{3,2}+3a_{2,1}a_{2,2})x_2$,\smallbreak\quad
$\nabla
R(\partial_{x_2},\partial_{x_3},\partial_{x_3},\partial_{x_2};\partial_{x_1})=
-\frac{2}{3}(-4+3a_{2,1}+3a_{3,1}+3a_{1,1}a_{1,2})x_1$,\smallbreak\quad
$\nabla
R(\partial_{x_2},\partial_{x_1},\partial_{x_1},\partial_{x_3};\partial_{x_2})
=(2-a_{1,1}-a_{1,2}+a_{2,1}-a_{2,2}$\par\qquad\qquad\qquad\qquad
$+a_{3,1}-a_{3,2}+a_{1,1}a_{1,2}-a_{2,1}a_{2,2}-a_{3,1}a_{3,2})x_3$,\smallbreak\quad
$\nabla
R(\partial_{x_2},\partial_{x_1},\partial_{x_1},\partial_{x_3};\partial_{x_3})=
(2-a_{1,1}-a_{1,2}+a_{2,1}-a_{2,2}$\par\qquad\qquad\qquad\qquad
$+a_{3,1}-a_{3,2}+a_{1,1}a_{1,2}-a_{2,1}a_{2,2}-a_{3,1}a_{3,2})x_2$,\smallbreak\quad
$\nabla
R(\partial_{x_1},\partial_{x_2},\partial_{x_2},\partial_{x_3};\partial_{x_1})
=(2-a_{1,1}+a_{1,2}-a_{2,1}-a_{2,2}$\par\qquad\qquad\qquad\qquad
$-a_{3,1}+a_{3,2}-a_{1,1}a_{1,2}+a_{2,1}a_{2,2}-a_{3,1}a_{3,2})x_3$,\smallbreak\quad
$\nabla
R(\partial_{x_1},\partial_{x_2},\partial_{x_2},\partial_{x_3};\partial_{x_3})=
(2-a_{1,1}+a_{1,2}-a_{2,1}-a_{2,2}$\par\qquad\qquad\qquad\qquad
$-a_{3,1}+a_{3,2}-a_{1,1}a_{1,2}+a_{2,1}a_{2,2}-a_{3,1}a_{3,2})x_1$,\smallbreak\quad
$\nabla
R(\partial_{x_1},\partial_{x_3},\partial_{x_3},\partial_{x_2};\partial_{x_1})=
(\frac{2}{3}+a_{1,1}-a_{1,2}-a_{2,1}+a_{2,2}$\par\qquad\qquad\qquad\qquad
$-a_{3,1}-a_{3,2}-a_{1,1}a_{1,2}-a_{2,1}a_{2,2}+a_{3,1}a_{3,2})x_2$,\smallbreak\quad
$\nabla
R(\partial_{x_1},\partial_{x_3},\partial_{x_3},\partial_{x_2};\partial_{x_2})=
(\frac{2}{3}+a_{1,1}-a_{1,2}-a_{2,1}+a_{2,2}$\par\qquad\qquad\qquad\qquad
$-a_{3,1}-a_{3,2}-a_{1,1}a_{1,2}-a_{2,1}a_{2,2}+a_{3,1}a_{3,2})x_1$.
\medbreak\noindent We set $\nabla R=0$ to obtain the desired equations of Theorem \ref{thm-1.9}; the first $3$ equations
generate the last $6$.\hfill\qedbox

\section*{Acknowledgments} Research of M. Brozos-V\'azquez
supported by the project BFM 2003-02949 Spain and by the Max Planck Institute for the Mathematical Sciences in Leipzig,
Germany. Research of P. Gilkey supported by the Max Planck Institute for the
Mathematical Sciences in Leipzig, Germany.
Research of S. Nik\v cevi\'c
supported by DAAD
 (Germany), TU Berlin, MM 1646 (Serbia) and by the
project 144032D (Serbia). This paper is
dedicated to the memory of our colleague Novica Bla{\v z}i{\'c} who
passed away Monday 10 October 2005.

\end{document}